# The Hidden Toll of COVID-19 on Opioid Mortality in Georgia: A Bayesian Excess Opioid Mortality Analysis

Cyen J Peterkin[a,*], Lance A Waller[a], Emily N Peterson[a]

[a]*Department of Biostatistics and Bioinformatics, Emory Rollins School of Public Health, Atlanta, GA, USA*

**Abstract**

COVID-19 has had a large scale negative impact on the health of opioid users exacerbating the health of an already vulnerable population. Critical information on the total impact of COVID-19 on opioid users is unknown due to a lack of comprehensive data on COVID-19 cases, inaccurate diagnostic coding, and lack of data coverage. To assess the impact of COVID-19 on small-area opioid mortality, we developed a Bayesian hierarchical excess opioid mortality modeling approach. We incorporate spatio-temporal autocorrelation structures to allow for sharing of information across small areas and time to reduce uncertainty in small area estimates. Excess mortality is defined as the difference between observed trends after a crisis and expected trends based on observed historical trends, which captures the total increase in observed mortality rates compared to what was expected prior to the crisis. We illustrate the application of our approach to assess excess opioid mortality risk estimates for 159 counties in GA. Using our proposed approach will help inform interventions in opioid-related public health responses, policies, and resource allocation. The application of this work also provides a general framework for improving the estimation and mapping of health indicators during crisis periods for the opioid user population.

*Keywords:* Bayesian Hierarchical models, COVID-19, Disease mapping, Excess Opioid Mortality

[*]Corresponding Author
  *Email address:* `cyen.peterkin@emory.edu` (Cyen J Peterkin)

# 1. Introduction

In March 2020, the World Health Organization declared COVID-19 a global pandemic, prompting the United States to implement a series of lockdowns to curb its spread. These measures disrupted the services and treatments available to individuals battling opioid use addiction, i.e., clinics providing methadone and other addiction treatment medications were forced to close, and access to addiction support groups became limited (Rossen et al., 2020; Leon et al., 2020). Direct and indirect impacts of the COVID-19 pandemic on opioid users resulted in a notable increase in opioid mortality rates. In Georgia, from 2019 to 2021, the total number of opioid-involved overdose deaths increased by 101% (Georgia Department of Public Health (GADPH), 2021). Critically, the total impact of COVID-19 on opioid users in Georgia, resulting from decreased utilization of routine health services, remains unknown. This information is a key ingredient in the resilience of the health system, essential for adequate allocation of resources to both 'crisis response activities' and 'core functions' (Sochas et al., 2017). However, monitoring opioid mortality during crisis periods is fraught with challenges that undermine a timely and accurate understanding of the crisis. A significant barrier is incomplete or inaccurate data, often due to underreporting, misclassification, and inconsistent use of ICD codes in cause-of-death certification (Stokes et al., 2021; Kline and Hepler, 2021; Hepler et al., 2021). Geographic and demographic disparities further exacerbate these gaps, leaving small areas, rural regions, and marginalized populations underrepresented in mortality statistics (Agarwal et al., 2002; Konstantinoudis et al., 2023). Development of robust and granular methodologies to monitor small-area geographical-temporal trends in opioid mortality during a crisis period has large-scale implications for informing public health policy, resource allocation, and public health response strategies (Jon Wakefield, 2007; Kline D and Hepler SA, 2021; Sochas et al., 2017). We present a Bayesian Excess Opioid Mortality Model with hierarchical spatio-temporal structures to estimate small-area opioid mortality rates during the COVID-19 crisis. Applied to Georgia's county-level data (2018–2022), the model offers robust and granular insights into the pandemic's impact on opioid users in the state.

Excess mortality has been used to assess the total impact of a crisis on public health outcomes when direct information is either sparse or unavailable, and accuracy is questionable. Excess mortality is defined as the difference between observed death counts for a given time period and the expected number of deaths based on historical time trends pre-crisis (Rossen et al., 2020; Woolf et al., 2020; Wang et al., 2022; Garfield, 2007; Checchi and Roberts, 2008; Stang et al., 2020; Blangiardo et al., 2020). Previous studies have primarily focused on assessing excess all-cause mortality due to COVID-19 at the national or sub-national levels. Woolf et al. (2020) analyzed state-specific excess deaths in the U.S. from COVID-19 and other causes during March–July 2020 using a Poisson regression model to estimate expected mortality. Similarly, Rossen et al. (2020) investigated excess deaths by age, race, and ethnicity in the U.S. for the period January 26–October 3, 2020, providing valuable demographic insights. Blangiardo et al. (2020) applied a spatio-temporal disease mapping approach to estimate weekly variations in excess mortality at the sub-national level in Italy, specifically evaluating excess mortality at the municipality level from January 1–April 28, 2020, and tracking its evolution over time. Vanella et al. (2021) used a principal component analysis to assess country-level weekkly mortality data by sex and age strata. Banerjee et al. (2020) estimates excess 1-year mortality associated with COVID-19 based on varying conditions and age. However, assessment of excess all-cause mortality at national and sub-national levels does not give a comprehensive understanding of the impact of COVID-19 on opioid users. To address this gap, our study focuses on excess opioid mortality, providing a more targeted understanding of the pandemic's impact on the particularly vulnerable population of opioid users. Small area estimation of excess opioid mortality at the county-level yields several statistical challenges: (1) an excess number of zeroes in extremely small counties, (2) very small population sizes compared to aggregates at the national and state levels, and (3) overdispersion in the data, driven by the combination of rare events and heterogeneity in underlying risk factors. Our approach not enhances the precision and robustness of opioid mortality estimates in the presence of sparse data, but also offers actionable insights for public health interventions aimed at mitigating the disproportionate burden of the opioid crisis.

Our proposed Bayesian Opioid Excess Mortality (BOEM) approach produces model-based county estimates of the expected number of opioid overdoses based on previous observed spatial-temporal trends pre-pandemic, which is crucial for evaluating excess opioid mortality rates (Rossen et al., 2020; Wang et al., 2022). The model framework employs a spatio-temporal disease mapping approach, which is commonly used to investigate granular geographical-temporal variations of opioid mortality burden (Kline D and Hepler SA, 2021; Hepler et al., 2021; Hepler SA et al., 2023; Kline D et al., 2023). The standard Bayesian disease mapping model provides the capability to borrow strength and share information across small areas thereby reducing high degrees of uncertainty associated with smaller population sizes (Waller, L.A. and Gotway, C.A., 2004; Banerjee et al., 2020; Andrew B. Lawson, 2013; Jon Wakefield, 2007). This is accomplished through spatial and temporal autocorrelation terms included within the hierarchical structure informing estimates of the unknown true risk of opioid-related death for a given county-time. Within our model framework, monthly opioid-related death counts between the years 2018-2019 are used to forecast county-month specific opioid mortality relative risk (RR) estimates for the years 2018-2022. The forecasted estimates capture the expected county-month specific opioid mortality trends in the absence of COVID-19. We derive estimates of excess mortality for 159 counties in Georgia by comparing the observed monthly number of opioid-related deaths to corresponding forecasted estimates representing the total impact of COVID-19 on opioid mortality rates. We illustrate model results for selected counties in Georgia to examine differing excess mortality trends across a wide variety of population sizes, and county characteristics. Using these findings, we aim to identify areas of greater need that require more intervention procedures, policies, and resource allocation to effectively address substance abuse, and counties suffering higher negative impact as a result of COVID-19.

This paper is organized as follows: Section 2.1 describes the data used to obtain excess opioid mortality estimates. Section 2.3 describes the data model assumed for the observed opioid death counts. Section 2.4 describes model assumption of the underlying process. Section 2.5 summarizes the process to obtain excess opioid mortality estimates and associated uncertainties. Lastly, Section 5 illustrates results across small, medium, and large population cases.

## 2. Methods

*2.1. Data*

Georgia Department of Public Health (GADPH) provided monthly opioid-related death counts for all 159 counties in Georgia for 2018 through 2022 that were aggregated to county-month specific totals (Georgia Department of Public Health, 2024). The U.S. Census Bureau (USCB) publishes annual county-level population estimates using a cohort component model based on the last available census with adjustments for births, deaths, and net-migrations (Population Estimation Program, U.S. Census Bureau, 2021). Given the USCB does not report monthly population estimates, we assume that the population size remains constant throughout each respective year (12 months). Dividing the county-month specific opioid-related death total by the corresponding population estimate gives the crude rate of opioid overdoses by county-month from 2018-2022. Figure 1 illustrates the reported total count of opioid-related deaths across Georgia from 2018-2022. The red line signifies the onset of COVID-19 (Jan. 2020), and each respective color represents year-specific data. Figure 1 shows a relatively stable trend in reported total number of opioid-related deaths during pre-pandemic years. In 2020, the number of deaths increased with a significant spike in March. Between 2020 and 2022, the number of deaths continues to rise over time, not reverting to pre-COVID-19 levels.

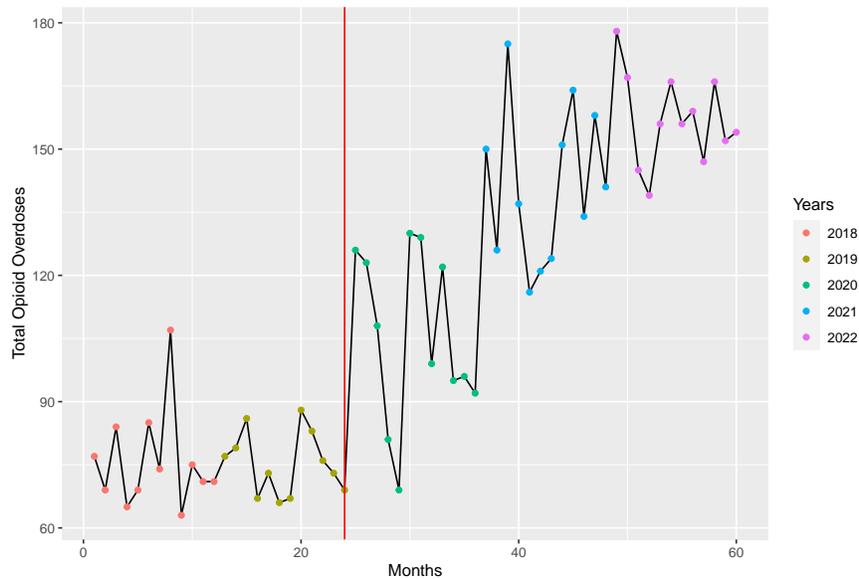

Figure 1: Monthly observed totals of opioid-related deaths across 159 counties in Georgia between the years 2018-2022. The red line denotes the defined start date of the COVID-19 pandemic (Jan. 2020). X-axis lists units in months ranging from 0 to 60, 0= Jan 2018, and 60 = Dec 2022. Colors denote years.

Figure 2 contains five maps of Georgia displaying the county-level crude rates of opioid related deaths (per 100,000 persons) by year. We note that higher rates occur in the years 2021 and 2022, a mass of higher rates is present mainly in the northwest region. The largest crude opioid mortality rate pre-pandemic occurred in Taliaferro County (Population: 1,537) in 2019. It is worth mentioning that Taliaferro County, being the least populated county in Georgia, results in highly noisy rate estimates. The largest crude opioid mortality rate during COVID-19 occurred in Wilkinson County (Population: 8,824) in 2021. These exploratory findings motivate our model assumptions.

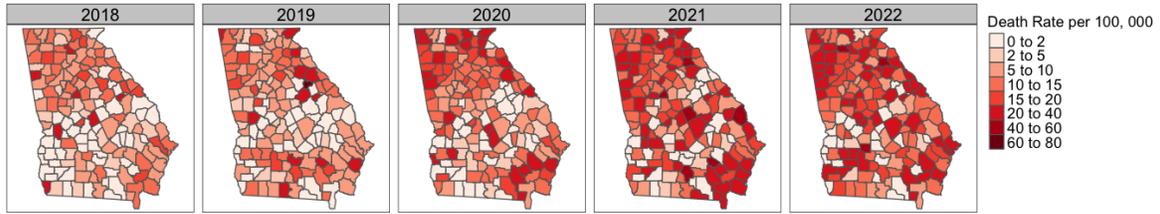

Figure 2: Mapped crude rates of opioid mortality (per 100,000) by year in Georgia between years 2018-2022. Color denotes the death rate per 100,000 which varies between 0 and 80.

## 2.2. Summary of model approach

Our Bayesian opioid excess mortality (BOEM) model estimates excess opioid mortality for 159 counties in Georgia to quantify the total (direct plus indirect) impact of COVID-19 on opioid users within the state. In our model approach we focus on estimating county-month specific opioid mortality rates, but note our approach could be applied to various spatial and temporal resolutions and socio-demographic sub-groups. The main features of the BOEM model are as follows:

1. The *data model* (defining the likelihood function) consists of modeling observed county-month specific opioid-related death counts for 2018-2019 $y_{im}$ using a Zero-Inflated Poisson assumption and is further detailed in Section 2.3.
2. The *process model* captures the latent true value of log-relative risk for county $i$ month $m$ denoted ($\theta_{im}$) which is modeled as a function of spatial and temporal random effects defining the assumed underlying process. This is further described in Section 2.4.
3. Estimates of latent and true excess opioid mortality denoted $\chi_{im}$ are derived using posterior predictive distribution estimates of expected opioid overdose counts, further detailed in Section 2.5.

Appendix A illustrates a graphical representation of the BOEM model set-up, refer to Appendix B for a summary of notation used.

## 2.3. Data model for observed opioid-mortality counts

In a standard disease mapping model, we commonly assume the relationship between the observed mortality counts and the true relative risk (i.e., the data model for observed cases) is given by a Poisson distribution which implicitly assumes cases in nearby areas are similar and the variance of response is equal to the mean not accounting for over-dispersion (Waller, L.A. and Gotway, C.A., 2004; Jon Wakefield, 2007; Andrew B. Lawson, 2013; Knorr-Held and Julian Besag, 1998; Parker PA et al., 2024). Given the nature of our observed county-level opioid mortality count data, which contains an excess of zero values particularly for smaller counties, we model observed counts using a Zero-Inflated Poisson (ZIP) model in Eq. 1, which accounts for over-dispersion due to the presence of excess zeros in the data (Agarwal et al., 2002; Ugarte et al., 2004; Feng, 2021). Let $y_{im}$ denote the observed counts for county $i$ and month $m$ (labeled $[im]$ for readability). The contribution of the Poisson likelihood is governed by the mixing parameter $\pi_{im}$, which represents the proportion of excess zeros. The county-month mixing parameter $\pi_{im}$ follows a Bernoulli distribution with global probability $\rho$, where $\rho$ is restricted to the interval [0,1]. It is important to note that $y_{im}$ represents the observed data from years 2018 and 2019, excluding data from 2020-2022. The full data model is given by:

$$f(y_{im}|\theta_{im}, X_{im}) \sim \begin{cases} \pi_{im} + (1-\pi_{im})e^{-(X_{im}\theta_{im})} & \text{if } y_{im} = 0 \\ (1-\pi_{im})\frac{\exp(-X_{im}\theta_{im})(X_{im}\theta_{im})^{y_{im}}}{y_{im}!} & \text{if } y_{im} > 0 \end{cases}, \quad \text{for } i = 1, ..., 159, , m = 1, ..., 24 \tag{1}$$

$$\pi_{im} \sim \text{Bern}(\rho)$$
$$\rho \sim \text{Unif}(0, 1)$$

The expectation $E(y_{im}|\theta_{im}, X_{im}, \pi_{im})$, shown in Eq. 2, is expressed as a function of the expected count $X_{im}$, the log relative risk $\theta_{im}$, and one minus the mixing parameter $\pi_{im}$. The offset $X_{im}$ is equal to the product of the reference rate $R$ and the population-at-risk $N_{im}$, both of which are considered fixed and known quantities. The reference rate is derived from a larger population base, resulting in substantially lower relative uncertainty compared to local estimates.

$$E(y_{im}|\theta_{im}, X_{im}) = \mu_{im} = (1-\pi_{im})X_{im}\theta_{im} \tag{2}$$
$$Var(y_{im}|\theta_{im}, X_{im}) = \mu_{im} + \frac{\pi_{im}}{1-\pi_{im}}\mu_{im}^2$$
$$X_{im} = R \cdot N_{im}$$
$$R = \frac{\sum_{im} y_{im}}{\sum_{im} N_{im}}$$

### 2.4. Process model for unobserved latent opioid mortality log-relative-risks

We model the latent log relative risk $\theta_{im}$ incorporating both spatially and temporally structured and unstructured random effects shown in Eq. 3. To incorporate spatial terms in our model, we consider the Besag-York-Mollie (BYM) Model (J Besag et al., 1991), which allows us to estimate the relative risk of death weighting trends in the neighboring counties. We denote $v_i$ to represent a spatially unstructured random effect term that is independent, identically, and normally distributed centered around zero, i.e., $v_i \sim N(0, \sigma_v^2)$). The spatially structured term, denoted $u_i$ is modeled assuming an Intrinsic Conditional Auto-Regressive (ICAR) prior, which assumes complete correlation between neighboring areas. The spatial covariance matrix, $W$ is written as a function of an $N \times N$ adjacency matrix where entries $\{i, i\}$ are zero and the off-diagonal elements are 1 if counties i and j are neighbors and 0 otherwise. $D$ is the $N \times N$ diagonal matrix where entries $\{i, i\}$ are the number of neighbors of county $i$ and the off-diagonal entries are 0. Lastly, $\tau$ denotes the smoothing parameter. It is important to note that our spatial parameters do not change over time. Correlated time trends are captured using a time structured random effect $\kappa_m$ modeled with a random walk order 1 which assumes a constant trend in forecasted estimates. To capture the non-separable process of space-time, we introduce $\omega_{im}$ which captures the interaction between terms $u$ and $\kappa$. As such the interaction term captures the deviations away from the separable space-time trend (Knorr-Held and Julian Besag, 1998). We model the interaction term as a Markov Random Field (MRF) where both first order and second order terms of neighbors are included, which assumes the temporal trend in county $i$ is similar to the average temporal trends of its neighbors (Knorr-Held, 2000). Using the complete process model, we obtain estimates of county-month specific log-transformed relative risks for years 2018-2022 including those months in which data has been excluded (during COVID months), i.e., 2020,...,2022.

$$
\begin{aligned}
log(\theta_{im}) &= \alpha + u_i + v_i + \kappa_m + \omega_{im}, \\
\alpha &\sim N(0, \sigma_\alpha^2) \\
v_i &\sim N(0, \sigma_v^2) \\
\boldsymbol{u} &\sim N(0, [\tau(D-W)]^{-1}) \\
\kappa_m &\sim N(\kappa_{m-1}, \sigma_\kappa^2) \\
\omega_{im} &= \omega_{i,m-1} + \sum_{j \in N(i)} \omega_{j,m} + \sum_{j \in N(i)} \omega_{j,m-1} + \varepsilon_{im} \\
\varepsilon_{im} &\sim N(0, \sigma_\omega^2)
\end{aligned}
\tag{3}
$$

for $i = 1,...,159, m = 1,...,60$
Intercept
Unstructured Spatial Noise
(ICAR) prior for spatial auto-correlation
Random Walk(1) for temporal autocorrelation
Space-time interaction (Type 4)
Space-time deviations

*2.5. Derivation of Excess mortality and associated uncertainties*

To derive excess mortality estimates, we use the posterior predictive distribution (PPD) to obtain predicted counts of opioid related death for each county-month (Konstantinoudis et al., 2023; Blangiardo et al., 2020). The PPD is represented by Eq. 4, where $p(\tilde{y}|y)$ represents the posterior distribution of expected counts $\tilde{y}$ given the observed data $y$, and the posterior estimates of the log-relative risk $p(\theta_{im}|y_{im})$.

$$
p(\tilde{y}|y) = \int_\theta p(\tilde{y}_{im}|\theta_{im}, y_{im}) p(\theta_{im}|y_{im}) d\theta_{im} \tag{4}
$$

Using posterior predictive sample estimates of opioid death counts, denoted $\tilde{y}^{(s)}$ for samples $s = 1,...,S$, excess mortality sample estimates $\chi_{im}^{(s)}$ are equal to the difference between the observed and predicted values shown in Eq. 5. Median and 95% uncertainty intervals are calculated by taking the median, and 95% quantile estimates of the sample values $\chi_{im}^{(s)}$ across all $S$ samples.

$$
\chi_{im}^{(s)} = \tilde{y}_{im}^{(s)} - y_{im} \text{ for } m = 1,...,60 \tag{5}
$$

## 3. Model Simulation

In this simulation exercise, we generate 100 synthetic data sets based on the proposed model structure given by Eqs. 1- 5 and compare the estimated parameters and excess mortality rates from the BOEM model outputs to the true values to assess model accuracy and robustness. Specifically, we generate county-month-specific opioid-related death counts using the predefined log-relative risk surface that incorporates spatial, temporal, and interaction effects given in Eq. 3 which assumes trends to remain consistent with pre-COVID levels. As such, the simulated data mimic real county-level time trends pre COVID-19 years. The BOEM model is then fit to the simulated data only for years 2018-2019. We calculate the difference between the model-based and true excess deaths ($error = \hat{\chi}_{im} - \chi_{im}$ for years with left-out data 2020-2022. Our ability to recover the true excess mortality estimates is evaluated using metrics such as bias, root mean squared error (RMSE), and credible interval coverage. This approach allows us to assess the model's accuracy, robustness, and sensitivity under varying levels of data sparsity, zero inflation, and noise. The full description of the simulation exercise and results can be found in Appendix C.

## 4. Computation

We extract PEP reported population estimates for 159 counties in Georgia, years 2018-2022, using the tidycensus package (Walker, K., 2020). For model processing and output, a Markov Chain Monte Carlo (MCMC) algorithm samples from the posterior distribution of the parameters via the software *Nimble* (de Valpine et al., 2018). Eight parallel chains were run with a total of 80,000 iterations in each chain. Of these, the first 40,000 iterations in each chain are discarded so the resulting chains contain 40,000 samples. Additionally, we thinned the samples to retain every 10th iteration after burn-in. Thus for each parameters there was a total 32,000 saved posterior samples. Standard diagnostic checks using traceplots were used to check for convergence (de Valpine et al., 2018; Gelman A. et al., 2013; Vehtari et al., 2017; Gabry et al., 2019).

## 5. Results

### 5.1. Simulation Results

Predictive performance results illustrate the BOEM model can robustly capture spatio-temporal trends in excess opioid mortality. In Table 1, we summarize the error measures between BOEM model-based estimates of the number of excess opioid deaths $\hat{\chi}_{i,m}$ and the true number of excess deaths $\tilde{\chi}_{i,m}$ generated from the procedure described in Section Appendix C. We derive out-of-sample predictive performance summary metrics for years 2020-2022 corresponding to years with removed data in which forward projections would be estimated based on past monthly data. The mean and median error were estimated at 0.07 and 0.052, respectively, corresponding to a relative error of 0.03% and 0.04%, respectively. This suggests robust BOEM model performance in estimating 2020-2022 forward projections. Coverage 95% prediction intervals were slightly conservative with inside coverage probability of 93%.

|  | Error | | | | Relative Error (%) | | | Inside 95% CI |
| ---: | ---: | ---: | ---: | ---: | ---: | ---: | ---: | ---: |
| N left-out | ME | MDE | MAE | MSE | MRE | MDRE | MARE | |
| 5,724 | 0.070 | 0.052 | 0.067 | 0.020 | 0.03 | 0.04 | 0.51 | 0.93 |

Table 1: Simulation results. The outcome measures are: Number of left-out observations (N left-out), mean error (ME), median error (MDE), median absolute error (MAE), mean squared error (MSE), mean relative error (MRE), median relative error (MDRE), median absolute relative error (MARE), as well as proportion left-out observations inside the respective 95% prediction intervals (PS) based on the training set.

### 5.2. Global Parameter Estimates

Global and hyper parameters consist of the global variance terms $\sigma^2_{()}$ and the global level $\alpha$. Table 2 shows posterior model-estimates and credible intervals for the global and hyper parameters. The global level $\alpha$, based on the global dataset of county opioid mortality counts was estimated to be -0.110, with a 95% credible interval (CI) spanning -0.237 to -0.011, which corresponds to the global average of the global dataset on the transformed scale. Posterior estimates of variance, and correlation components are also given in Table 2.

| Parameters | Mean | SD | 95% CI Lower Bound | Median | 95% CI Upper Bound |
|---|---|---|---|---|---|
| $\alpha$ | -0.114 | 0.057 | -0.237 | -0.110 | -0.011 |
| $\sigma_\alpha$ | 0.195 | 0.032 | 0.071 | 0.213 | 0.364 |
| $\sigma_v$ | 0.176 | 0.061 | 0.068 | 0.174 | 0.293 |
| $\sigma_u = 1/\tau^2$ | 0.445 | 0.131 | 0.179 | 0.446 | 0.699 |
| $\sigma_\kappa$ | 0.014 | 0.012 | 0.001 | 0.011 | 0.047 |
| $\sigma_\xi$ | 0.058 | 0.016 | 0.030 | 0.060 | 0.086 |
| $\sigma_\varepsilon$ | 0.022 | 0.009 | 0.012 | 0.025 | 0.031 |

Table 2: Global and hyper parameters included in the random effect terms of the (BOEM) model.

### 5.3. State trends of excess mortality in Georgia

Figure 3 presents five maps of Georgia showing posterior median excess opioid deaths per county by year. Compared to 2018-2019, the 2021-2022 maps reveal increased opioid deaths in Metro Atlanta, particularly in DeKalb, Fulton, and Gwinnett counties, following COVID-19. In 2018-2019, excess deaths ranged from -1 to -3, indicating observed deaths were lower than predicted. By 2022, DeKalb and Fulton had 10 excess deaths each, while Gwinnett saw similar increases in 2020 and 2021. These findings highlight the stark shift in opioid deaths after the onset of COVID-19 in these high-risk areas.

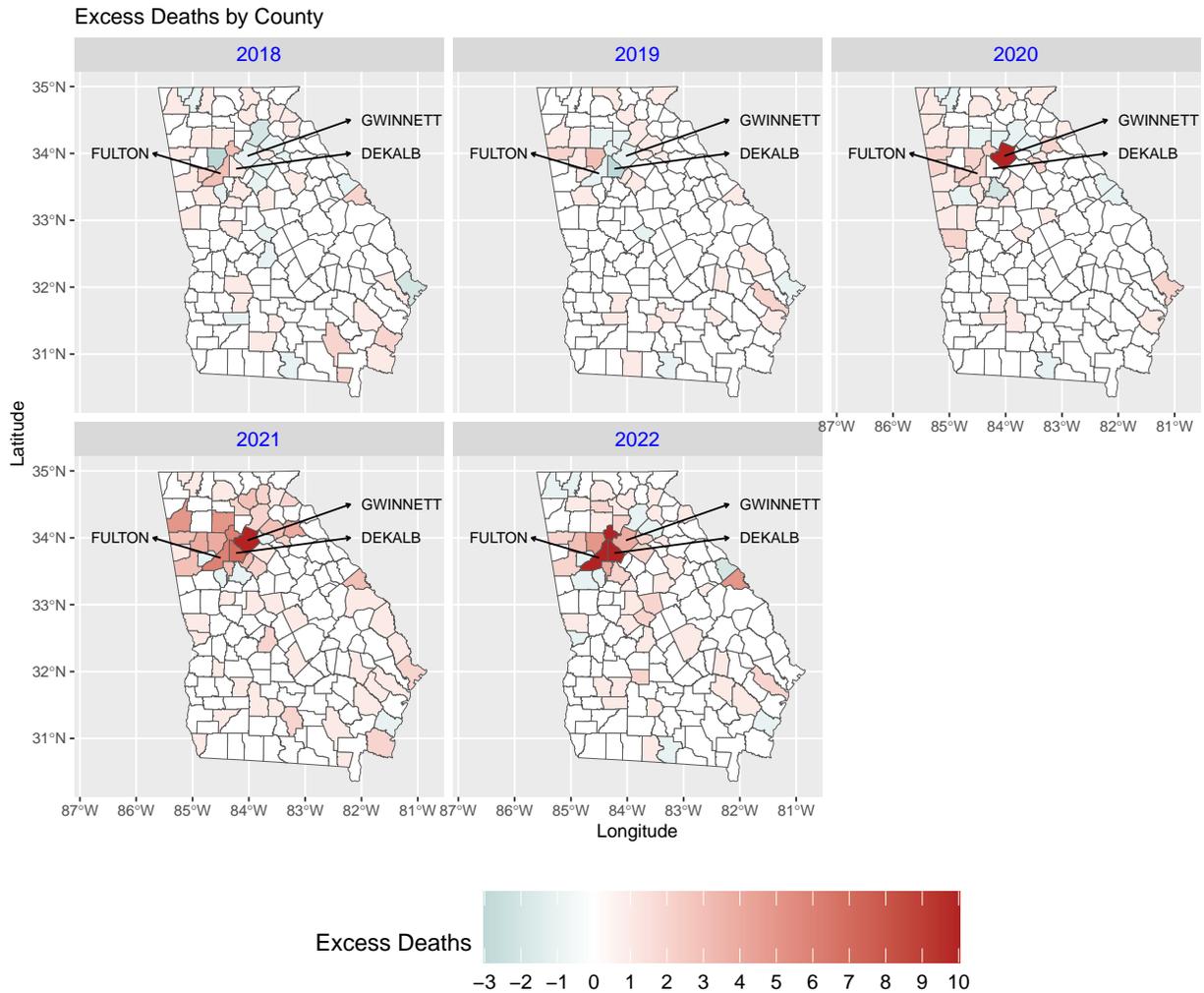

Figure 3: Mapped excess opioid deaths broken down by year in Georgia between years 2018-2022. The divergent color scale highlights areas with varying levels of excess mortality.

Figure 4 presents five maps of Georgia displaying posterior median excess opioid mortality rates (EMRs) per 100,000 persons by year. EMRs account for the relative size of the population-at-risk by dividing excess deaths by the population size. Pre-COVID years (2018-2019) show low or negative EMRs in many counties, indicating fewer deaths than expected (e.g., Appling County had an EMR of 0). Post-COVID (2020-2022) reveals significant increases in EMRs, especially in 2021-2022. For instance, Baker County reached a median EMR of 35 (95% credible bounds: 30–40) in 2021, and Berrien County recorded a median EMR of 22 (95% credible bounds: 18–27), reflecting extreme localized crises.

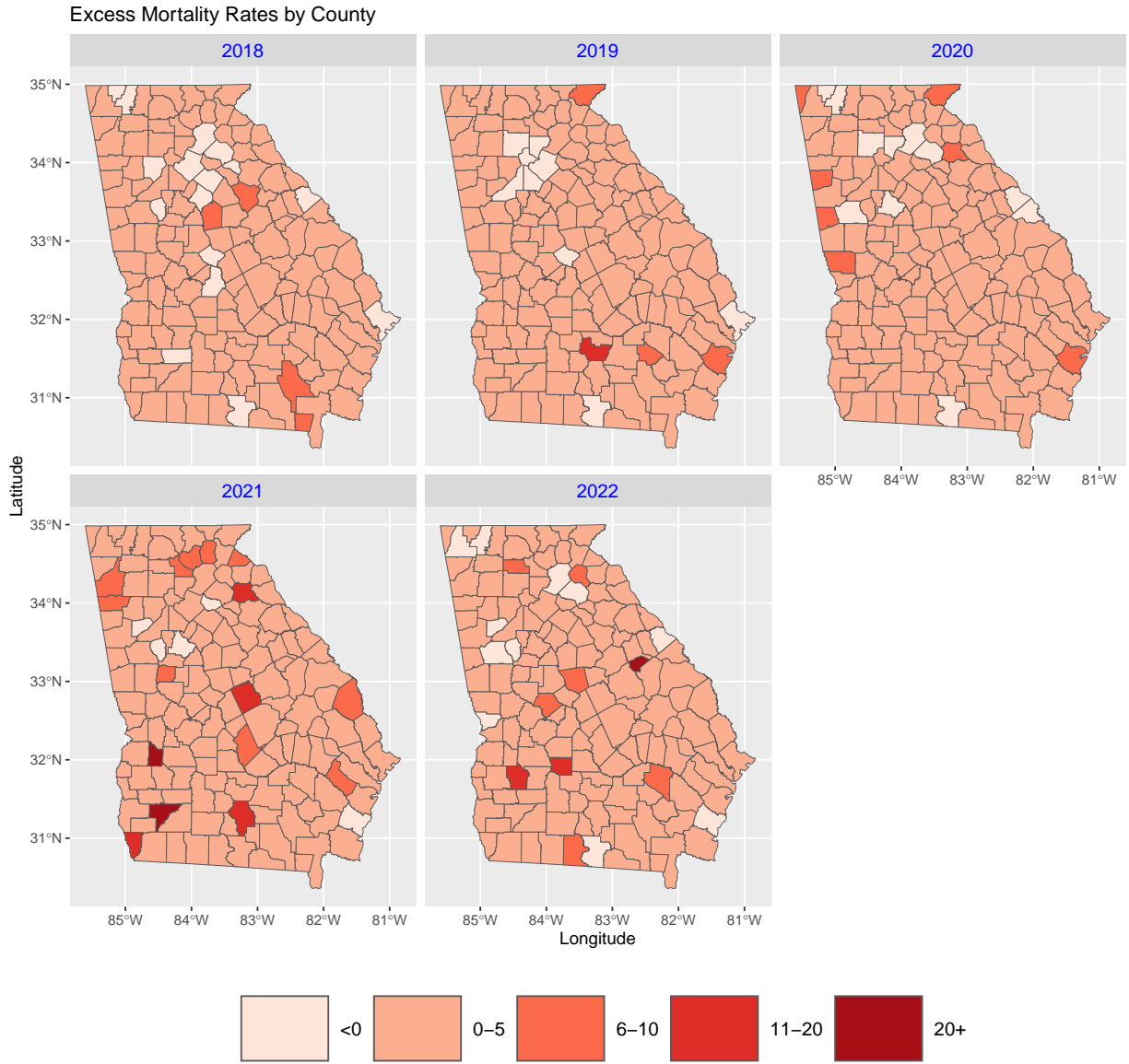

Figure 4: Mapped excess opioid deaths broken down by year in Georgia between years 2018-2022. The divergent color scale highlights areas with varying levels of excess mortality.

## 5.4. Large County Population Cases

Figure 5 illustrates excess opioid mortality trends for selected largely populated counties within Georgia including DeKalb, Fulton, and Gwinnett County. Each respective county plot represents the predicted (blue), observed (red), and excess (orange) rates of opioid mortality rate (per 100,000) by month for years 2018-2022.

In DeKalb County (Population: 754,906 - 764,420), the observed monthly opioid-related death rates per 100,000 range between 0 and 1.06 for years 2018-2019. The monthly opioid death rate rose to 1.97 in December 2022. The lowest predicted rate, which was 0.53 (95% CI: 0.27, 0.79) in May 2018, increased to 0.66 (95% CI: 0.27, 1.70) in 2022. These predictions for 2020-2022

were influenced by the observed trend from 2018-2019. The 95% uncertainty bounds around the predictive estimates widen as the predictions extend further from December 2019. This increase in uncertainty is expected because the predictions are based on data that becomes progressively older, making future estimates more uncertain. DeKalb County experienced its highest excess death rate of 1.45 (95% CI: 0.66, 1.84) in May 2021 and its lowest excess death rate of -0.40 (95% CI: -0.66, -0.27) in February 2018. There is a notable trend in excess opioid mortality rates, which are consistently banded around 0 pre-COVID and consistently larger than 0 during the COVID period.

In Fulton County (Population: 1,050,131 - 1,069,370), the observed monthly opioid-related death rates per 100,000 range between 0.38 and 1.22 for years 2018-2019. The monthly opioid death rate increased to 2.16 in May 2021 showing a stark increase in the number of opioid deaths in 2021-2022. The lowest predicted rate 0.47 (95% CI: 0.28, 0.75) occurred in July 2019, which increased to 0.75 (95% CI: 0.37, 1.95) in 2022 based on 2018-2019 observed trends. The estimated highest excess death rate of 1.45 (95% CI: 0.66, 1.84) occurred in May 2021 and its lowest excess death rate of -0.29 (95% CI: -0.57, 0) in July 2018.

In Gwinnett County (Population: 927,337 - 975,353), the observed monthly opioid-related death rates per 100,000 range between 0 and 1.4 for years 2018-2019. The monthly opioid death rate increased to 1.84 in October 2022 showing a similar increasing trend in opioid deaths in years 2021-2022. The lowest predicted rate 0.43 (95% CI: 0.32, 0.65) occurred in May 2018, which increased to 0.62 (95% CI: 0.20, 1.53) in 2022. The estimated highest excess death rate of 1.23 (95% CI: 0.41, 1.64) occurred in October 2022 and its lowest excess death rate of -0.43 (95% CI: -0.65, -0.32) in May 2018.

A comparison of excess mortality rates across the selected large counties reveals a consistent trend of rising opioid mortality rates following the onset of the COVID-19 pandemic. Each county experienced a peak in excess deaths during the pandemic period, highlighting its significant impact. Notably, negative EMRs indicate periods where the predicted number of deaths was higher than the observed value.

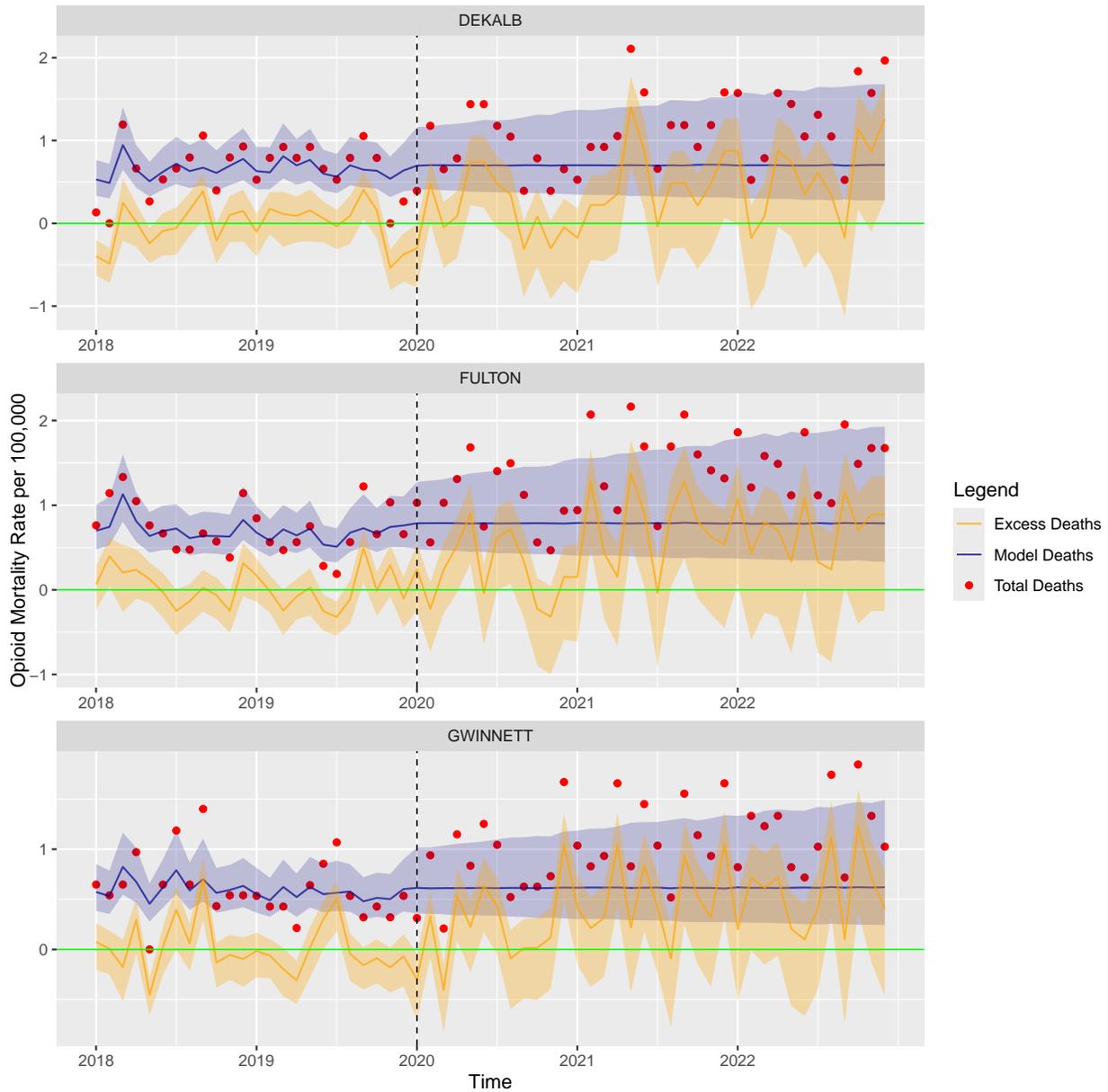

Figure 5: Monthly predicted, observed, and excess (and associated 95% confidence intervals) rates of opioid mortality (per 100,000) across DeKalb, Fulton, and Gwinnett County. The dashed line denotes the defined start of the COVID-19 pandemic in the model. Line colors distinguish the data source, red = observed deaths, blue = predicted deaths, orange = excess deaths, and green = horizontal line at 0.

## 5.5. *Moderate County Population Cases*

Figure 6 illustrates selected cases of moderate (medium) populated counties within Georgia including Bartow, Cherokee, and Clayton County.

In Bartow County (Population: 106,378 - 112,816), the observed monthly opioid-related death rates range between 0 and 4.64 for 2018–2019. This increases to 8.11 in February 2021. The larger variability in observed opioid death rates compared to larger counties illustrates the impact of smaller population sizes on observed fluctuations. Smaller populations suffer from stochastic variability in which a small increase in deaths results in substantial changes in observed rates. The lowest predicted rate 0.79 (95% CI: 0.45, 1.29) occurred in February 2018. The predicted rate increased to 1.25 (95% CI: 0.50, 3.27) in 2022. Bartow County experienced its highest excess death rate of 6.85 (5.34, 7.52) in February 2021 and its lowest excess death rate of -1.25 (95% CI: -2.29, -0.68) in January 2020.

In Cherokee County (Population: 253,914 - 281,278), the observed monthly opioid-related death rates per 100,000 range between 0 and 0.77 for years 2018-2019. The monthly opioid death rate increased to 2.13 in May 2022. The lowest predicted rate 0.66 (95% CI: 0.39, 1.045) occurred in July 2019, which increased to 0.90 (95% CI: 0.40, 2.10) in 2022. The estimated highest excess death rate of 2.33 (95% CI: 1.32, 2.88) occurred in March 2018 and its lowest excess death rate of -0.89 (95% CI: -1.56, -0.49) in January 2020. The spike in excess opioid deaths in March 2018 resulted from an increase in opioid deaths from 3 in February to 10 in March, illustrating the sensitivity of small and moderate-sized counties to changes in monthly death counts.

In Clayton County (Population: 289,197 - 297,623), the observed monthly opioid-related death rates range between 0 and 1.03 for 2018–2019. The observed monthly death rate increased to 2.36 in April 2022. The predicted opioid mortality rate stayed consistently around 0.50 (95% CI: 0.20, 0.70). The estimated highest excess death rate of 1.85 (95% CI: 1.12, 2.15) occurred in April 2022 and its lowest excess rate of -0.51 (-1.28, -0.20) occurred in August 2022. Excess opioid death rates in Clayton County exhibited erratic fluctuations over time, with notable spikes in rates for years 2020-2022 showing a substantial increase in deaths during the COVID-19 pandemic.

In all three counties, excess mortality rates increased following the onset of COVID-19, reflecting trends observed in larger counties, as shown in Figure 5. However, smaller and moderate-sized counties experience greater variability in both observed and predicted estimates due to their smaller population sizes. The sharp increase in opioid death rates in Cherokee County in March 2018 deviates from the overall trend in observed counts, highlighting how minor increases in deaths can lead to substantial fluctuations in death rates for smaller counties. These findings underscore the advantage of the BOEM approach, which leverages information sharing across counties to reduce uncertainty in estimates for smaller populations and introduces a level of smoothing to mitigate erratic fluctuations.

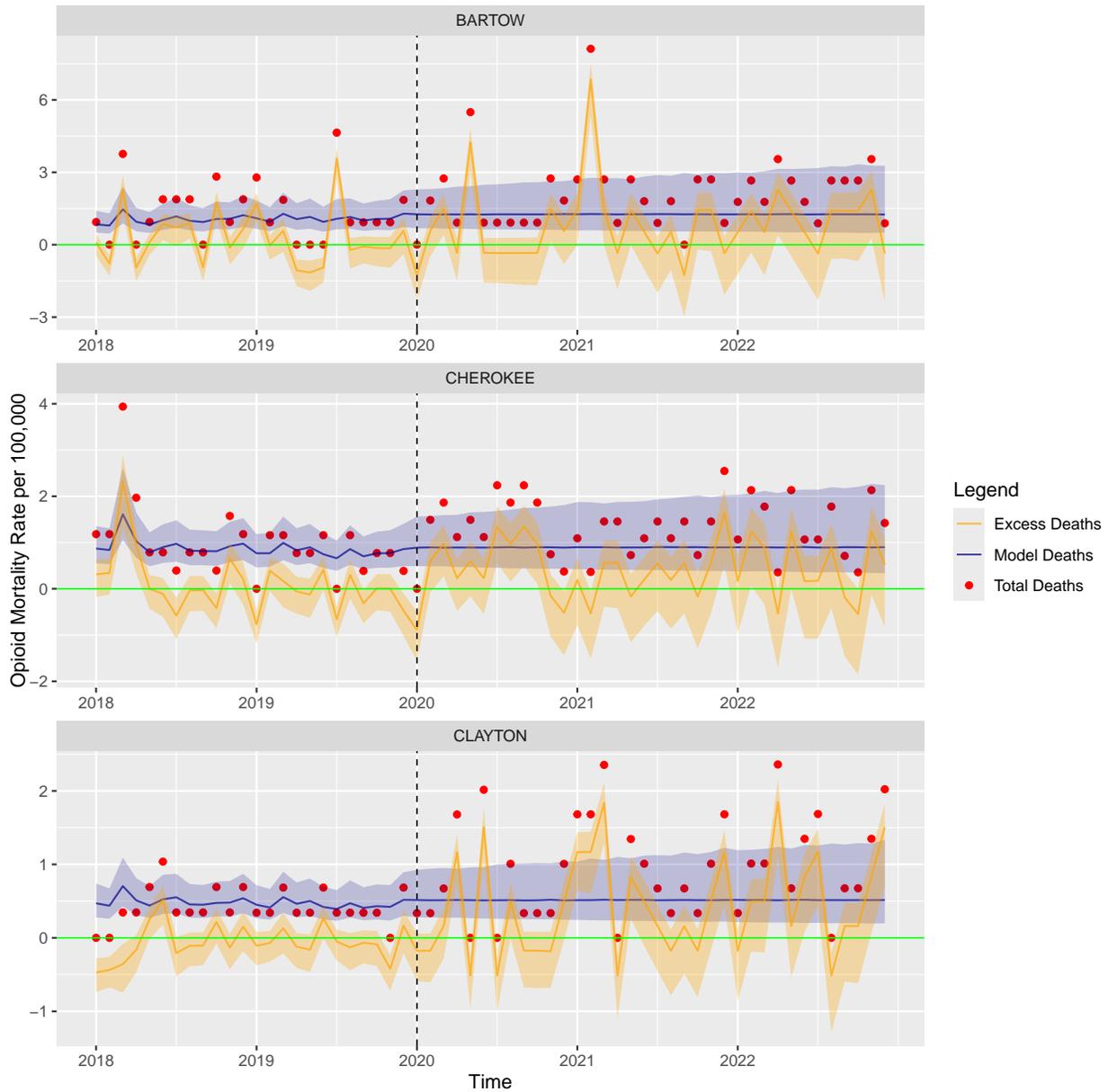

Figure 6: Monthly predicted, observed, and excess (and associated 95% confidence intervals) rates of opioid mortality (per 100,000) across Bartow, Cherokee, and Clayton County. The dashed line denotes the defined start of the COVID-19 pandemic in the model. Line colors distinguish the data source, red = observed deaths, blue = predicted deaths, orange = excess deaths, and green = horizontal line at 0.

## 5.6. Small County Population Cases

Figure 7 contains selected cases of small populated counties within Georgia including Toombs, Upson, and Walker County.

In Toombs County (Population: 26,830 - 27,081), the observed monthly opioid-related death rates stays at 0 for years 2018-2019. The observed trend erratically increases for years 2020-2022, with the highest observed death rate at 7.45 in September 2022. The predicted monthly opioid-related death rate maintains around 0.60 (95% CI: 0.20, 1.50). Toombs County experienced its highest excess death rate of 6.86 (5.75, 7.26) in September 2022 and its lowest excess death rate of -0.76 (95% CI:-1.32, -0.41) in March 2018. Toombs County serves as an example of a county with an excess number of zero deaths across the 24 months included in the data. Consequently, the predicted trend under the Zero-Inflated Poisson likelihood model is heavily influenced by the high frequency of zero counts. As a result, the predicted trend remains close to 0 across the complete time series.

Similar to Toombs County, Upson County (Population: 26,185 - 28,086) is another county with an abundance of observed zeroes across the complete time series. In contrast to Toombs County, there are spikes in the observed death rate both pre and post the onset of COVID-19. The largest observed death rate 3.18 occurred in August 2018. The predicted opioid mortality rate ranged from 0.49 (95% CI: 0.26, 1.43) in February 2018 to 0.62 (95% CI: 0.18, 1.76) in December 2022 showing a very slight increased in predicted trends from 2018 to 2022. The largest estimated EMR 3.28 (95% CI: 2.81, 3.55) occurred in August 2018.

In comparison Walker County (Population: 67,742 - 69,761) experienced fluctuating deaths rates from 2018 to 2022. The largest observed death rate 4.35 occurred in May 2022. The predicted opioid death rates ranged from 0.62 (95% CI: 0.36, 1.05) in January 2018 to 0.92 (95% CI: 0.33, 2.65) in January 2022. The highest predicted rate 1.03 (95% CI: 0.61, 1.74) occurred in March 2018. The largest estimated EMR of 3.46 (95% CI: 2.18, 3.98) occurred in May 2021. The lowest estimated EMR -0.92 (95% CI: -1.95, -0.44) occurred in June 2020. Unlike Toombs and Upson counties, which exhibit consistently low or zero opioid-related deaths, Walker County experienced more pronounced fluctuations in both observed and predicted rates, highlighting the greater variability and sensitivity to changes in opioid-related mortality in counties with small population sizes.

Figure 7 illustrates the stark and rapid changes in observed death rates that often occur due to singular deaths during particular times. Additionally, in these smaller counties, there are zero deaths for the majority of months. As such, excess mortality estimates obtained from the BOEM model capture the negative excess mortality rates in periods where there are zero deaths, and the increase in excess mortality in months where there are deaths present.

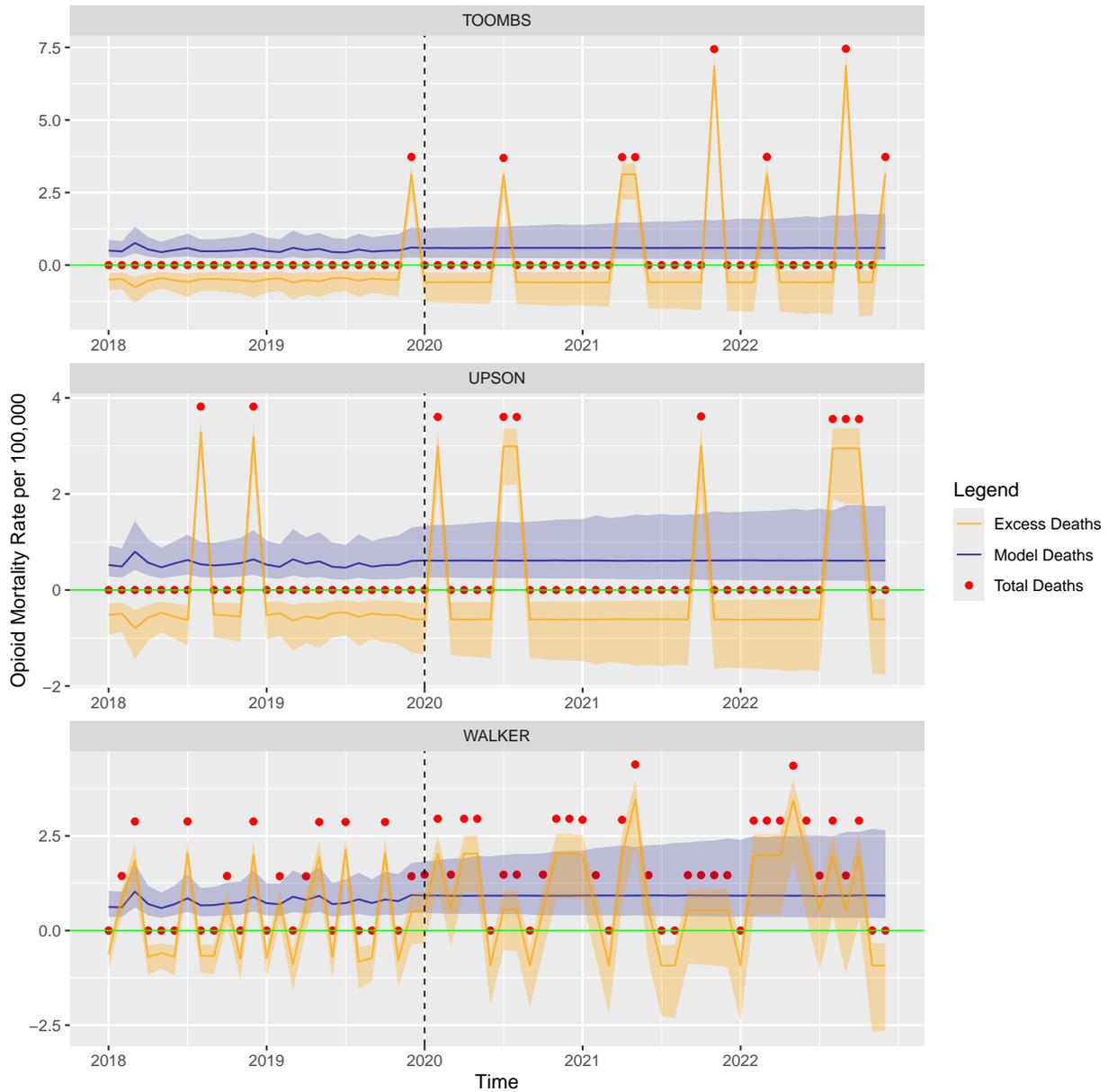

Figure 7: Monthly predicted, observed, and excess (and associated 95% confidence intervals) rates of opioid mortality (per 100,000) across Toombs, Upson, and Walker County. The dashed line denotes the defined start of the COVID-19 pandemic in the model. Line colors distinguish the data source, red = observed deaths, blue = predicted deaths, orange = excess deaths, and green = horizontal line at 0.

## 6. Discussion

Evaluating the impact of COVID-19 on opioid drug use is difficult due to the lack of real time cause-specific death data and inaccurate record keeping of COVID-19 cases. To address this, we have presented a Bayesian hierarchical excess opioid mortality model in which: (1) we explore the spatio-temporal variations in opioid deaths for 159 counties in Georgia, and (2) we assesses excess opioid mortality pre- and post- the onset of the COVID-19 pandemic. The general findings suggest that COVID-19 did act as a catalyst in excess opioid mortality. Pre-COVID-19, Georgia's monthly excess opioid related deaths did not exceed fifty. Post-COVID-19, Georgia's monthly excess opioid related deaths have spiked over fifty and forecasts to continue to increase. The county population size cases display key features of our model's performance and reactions to death counts, which vary based on the county's population size. In large counties, findings showed that they experienced an increase in their excess death rates post-COVID-19. In moderate size counties, larger variability is observed in their excess death rates due to the impact of deaths on their smaller population size. In smaller counties, the erratic nature of opioid mortality deaths illustrated the stark change in death rates due to a small number of individual events. Our findings guide understanding of the total impact of COVID-19 on opioid users within each county.

The contributions of this work to statistical modeling and excess opioid mortality are two-fold. Firstly, we expanded upon existing work that was done to evaluate excess mortality geographically as a result of COVID-19 to produce model-based excess mortality estimates at granular levels. Secondly, our developed methodology can be applied to assist in small area excess mortality estimation in other applications including other epidemics/pandemics, natural disasters, law regulations, mental health disorders, and more.

We also note the limitations of our study. Firstly, we observed a high number of zero values, which resulted in stationary estimates over time. In an attempt to address this issue, we replaced the use of a Poisson distribution in the standard disease mapping approach with a Zero-Inflated Poisson distribution. Secondly, in general, there is variability in forecasting and predicting time trends based on limited data.

**Funding** This work was supported by the Injury Prevention Research Center at Emory University Rollins School of Public Health under Grant Number R49CE003072.

# Appendix A. Graphical representation of the BOEM model

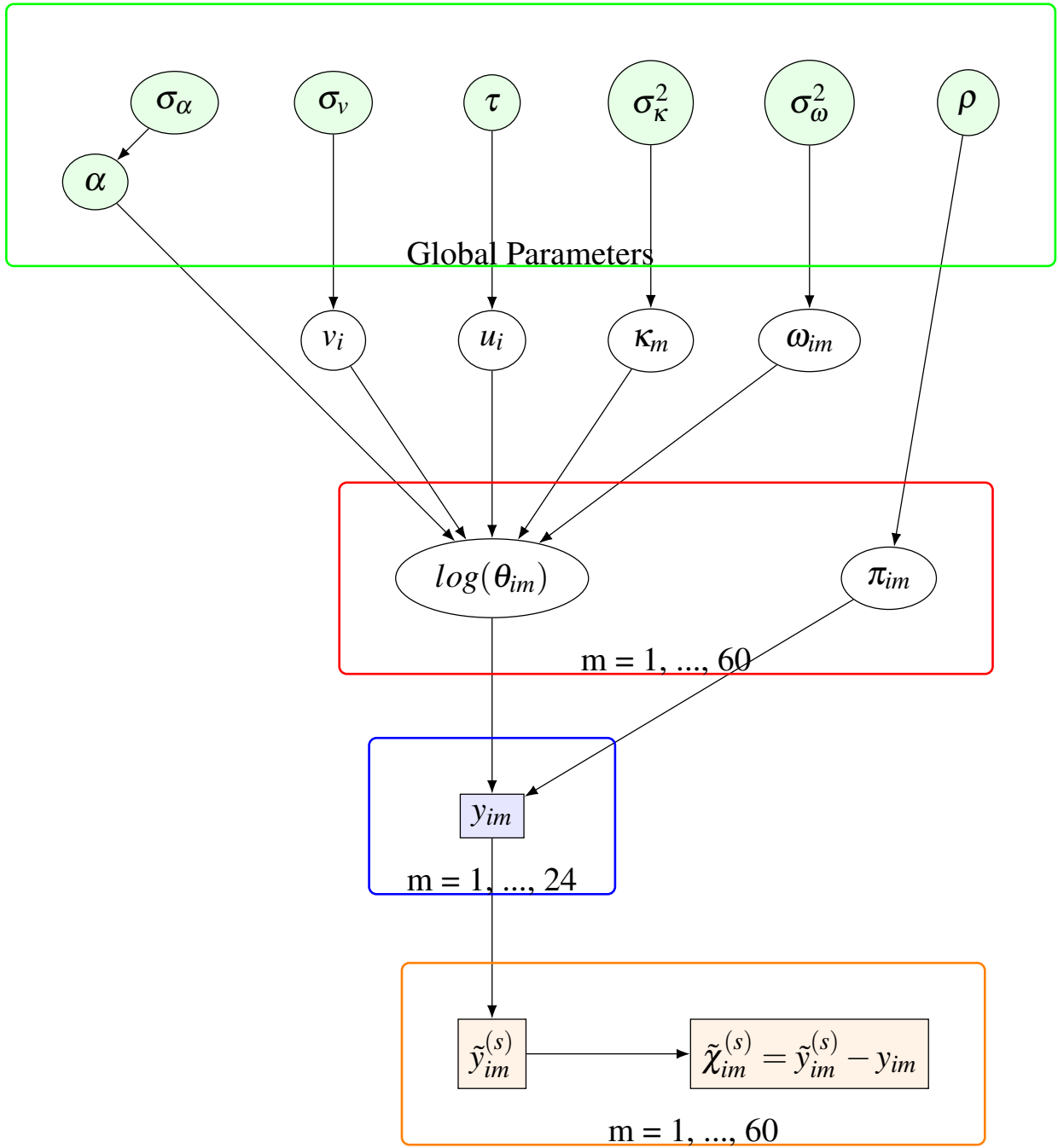

Figure A.8: Directed graphical representation of the BOEM hierarchical model. Shaded rectangles denote observed data quantities, and circles denote latent variables (shaded circles for global hyper-parameters). Solid arrows denote stochastic dependency. Boxes group quantities by indices, i.e., (1) Bottom box contains observed population data, stratified by county-month for years 2018-2019, (2) Middle box contains estimated parameters stratified by county-month, for years 2018-2022, (3) Top box contains global parameters. Subscripts refer to county $i$, month $m$.

## Appendix B. Notation Table

| Parameter Notation | Description |
|---|---|
| **Data Quantities** | |
| $R$ | Reference risk |
| $N_{im}$ | Population for county $i$, month $m$. |
| $y_{im}$ | Opioid mortality counts for county $i$, month $m$. |
| $\theta_{im}$ | Log relative risk for county $i$, month $m$. |
| $X_{im}$ | Population at risk for county $i$, month $m$. |
| $\pi_{im}$ | Poisson likelihood mixing parameter for county $i$, month $m$. |
| **Estimates Quantities** | |
| $\alpha$ | Overall Intercept. |
| $\sigma_\alpha^2$ | Variance of the intercept. |
| $v_i$ | Spatially unstructured random effect term. |
| $\sigma_v^2$ | Variance of the spatially unstructured random effect. |
| $u$ | Spatial auto-correlation (ICAR prior). |
| $\tau$ | Smoothing parameter of the ICAR prior. |
| $D$ | Diagonal matrix containing the number of neighbors of each area on the diagonal. |
| $W$ | Adjacency matrix containing 1 for neighboring & 0 for non-neighboring counties. |
| $\kappa_m$ | Temporal autocorrelation (RW(1)). |
| $\sigma_\kappa^2$ | Variance of the temporal structured random effect. |
| $\omega_{im}$ | Space-time interaction term (Type 4 spatio-temporal interaction). |
| $\sigma_\omega^2$ | Variance of the interaction between temporal and spatial effects. |
| $\varepsilon_{im}$ | Space-time deviation term modeled as independent and identically distributed normal distribution. |
| $\sigma_\omega^2$ | Space-time deviation variance term. |

## Appendix C. Simulation Description

*Data Generation.* We first generate log-relative risks for each county-month from the true latent process described in Eq. 3, incorporating fixed global parameter values. As such, we generate true log-relative-risks assuming no change due to the onset of COVID-19. The generated log-relative-risks capture the "truth" under the assumption that trends in 2018-2019 continue for years 2020-2022. The true process model includes the global level $\alpha$, the temporal term, $\tilde{\kappa}_m$, which follows a random walk (RW(1)) structure, the spatial structured and unstructured terms $\tilde{u}_i + \tilde{v}_i$, and the space-time interaction term, $\tilde{\omega}_{i,m}$, which is modeled as the product of spatial and temporal main effects. Additionally, a random noise term, $\tilde{\varepsilon}_{i,m}$, is included. The true underlying process model is defined in Eq. C.1, which is used to derive true annual county-specific opioid mortality relative-risk $\tilde{\theta}_{i,m}$ for all months from 2018 to 2022. Additionally, we derive the simulated county-month values of the mixing parameter $\tilde{\pi}_{i,m}$ assuming a Bernoulli distribution with a fixed global value for the probability $\rho$. Generated data quantities are denoted by $\tilde{x}$. A summary of the data generation process is provided below:

$$log(\tilde{\theta}_{i,m}) = \alpha + \tilde{u}_i + \tilde{v}_i + \tilde{\kappa}_m + \tilde{\omega}_{i,m} + \tilde{\varepsilon}_{i,m} \tag{C.1}$$
$$\alpha = 0.1 \text{ the global mean across the complete dataset}$$
$$\tilde{u}_{1:C} \sim ICAR(\tau = 1)$$
$$\tilde{v}_i \sim N(0, 0.1^2)$$
$$\tilde{\kappa}_m \sim RW(1)$$
$$\tilde{\omega}_{i,m} \sim N(0, 0.5^2)$$
$$\tilde{\pi}_{i,m} \sim Bern(\rho = 0.2)$$

We generate the "true" opioid mortality counts for county $i$ month $m$ $\tilde{y}_{i,m}$ assuming the zero-inflated Poisson data generating model given in Eq. 1. The expectation is equal to the product of the log-relative-risk $\tilde{\theta}_{i,m}$, the observed offset $X_{i,m} = R \cdot N_{i,m}$, and the generated mixing parameter $\tilde{\pi}_{i,m}$ obtained from Eq. C.1, i.e., $\tilde{\theta}_{i,m} \cdot X_{i,m} \cdot \tilde{\pi}_{i,m}$. From the generated true opioid mortality counts $\tilde{y}_{i,m}$, we can derive the true number of excess opioid deaths defined as the difference between the observed number of deaths after the onset of Covid-19 $y_{i,m}$ and the true number of deaths generated under the assumed true data generating assumption without a change due to COVID-19 $\tilde{y}_{i,m}$, i.e., $\tilde{\chi}_{i,m} = y_{i,m} - \tilde{y}_{i,m}$ which captures the true number of excess opioid deaths for a given county-month.

*Model Fitting.* In our simulation exercise, we generate 100 datasets of simulated county-month specific relative risks and corresponding opioid mortality counts using the model described in Eq. C.1, which assumes temporal trends differ across counties, but trends in adjacent counties are more similar. We restrict the generated population counts to be available only for years 2018-2019. We fit the BOEM model summarized in Sections 2.1-2.4 to the 100 generated datasets of true opioid mortality counts (2018-2019) to obtain model-based errors defined as the difference between the BOEM estimated excess mortality rate $\hat{\chi}_{i,m}$ and the true excess mortality rate $\tilde{\chi}_{i,m}$. We assess model performance using summary metrics of mean error, median error, median absolute error, mean squared error, relative errors, and 95% coverage intervals. The summary of the model fitting and model performance assessment is given in Figure C.9.

**Calculation of outcome measures in the validation exercise**

1. Model fitting: Fit the BOEM model to the training data ($m = 1, ..., 24$) for years 2018-2019 and obtain posterior samples $\theta_{i,m}^{(s)}$ for county-months with left-out data in the test set ($m = 25, ..., 60$) for years 2020-2022.
2. Derive the estimated number of excess deaths for months ($m = 25, ..., 60$) defined as the median of the posterior predictive distribution (PPD) across posterior samples $s = 1, ..., S$.

$$\hat{y}_{i,m} = median(\hat{y}_{i,m}^{(s)}) \text{ for } m = 25, ..., 60$$

3. Calculate the estimated number of excess opioid deaths defined as the difference between the observed number of opioid deaths and the posterior median of the PPD.

$$\hat{\chi}_{i,m} = y_{i,m} - \hat{y}_{i,m} \text{ for } m = 25, ..., 60$$

4. Error calculation: Calculate the difference between estimated number of excess deaths and the true number of excess deaths based on simulated data.

$$error_{i,m} = \hat{\chi}_{i,m} - \tilde{\chi}_{i,m} \text{ for } m = 25, ..., 60$$

Various summaries of the errors are reported.

5. Calibration: Calculate the proportion of observed excess deaths counts above and below their respective 95% prediction interval.

Figure C.9: Overview of calculation of errors and coverage of prediction intervals in out-of-sample validation exercises.